\documentclass{mysubmit}
\newcommand{\pW}{\mathop{\mathit{pW}}}
\newcommand{\pP}{\mathop{\mathit{pP}}}

\newcounter{myFigure}
\newenvironment{myFigure}{\begin{center}}{\end{center}}
\renewcommand{\caption}[1]{\refstepcounter{myFigure}\textsf{Figure~\arabic{myFigure}. }\emph{#1}}

\newtheorem{claim}{Claim}

\newcommand{\calG}{\mathcal{G}}

\begin{document}

\begin{center} 
{\Large \bfseries
Coloring Graphs to Produce Properly Colored Walks} \\[6mm]
\begin{tabular}{c}
{\large Robert Melville  and
Wayne Goddard} \\[3mm]
Dept of Mathematical Sciences \\
Clemson University \\
\textsf{goddard@clemson.edu}
\end{tabular}
\end{center}

\begin{abstract}
For a connected graph, we define
the proper-walk connection number as the 
minimum number of colors needed to color the edges of a graph so that 
there is a walk between every pair of vertices without 
two consecutive edges having the same color. 
We show that the proper-walk 
connection number is at most three for all cyclic graphs,
and at most two for bridgeless graphs.
We also characterize the bipartite graphs that 
have proper-walk connection number equal to two, and show that this 
characterization also holds for the analogous problem
where one is restricted to properly colored paths.
\end{abstract}

\section{Introduction}

We consider the problem of coloring the edges of a graph so that 
it is possible to get between every pair of vertices without 
two consecutive edges having the same color. 
Obviously, this can be achieved by
giving every edge a different color, and indeed by
any proper coloring of the edges. So the real question is 
what is the minimum number of colors one needs.

Borozan et al.~\cite{BFGMMMT12proper} introduced this problem
for paths. In particular, for a connected graph $G$, they defined  
the \emph{proper connection number} as the 
minimum number of colors that one needs so that there is 
a properly colored path between every two vertices. For example,
they showed that the parameter is at most $3$ for any block. 
Also, if a graph has a Hamiltonian path,
then the parameter is at most $2$~\cite{ALLZ16proper},
and thus almost surely this holds for a random graph~\cite{GLQ16random}.
For a recent survey, see~\cite{LM15survey}.

We consider here the analogous concept for walks.
For a connected graph $G$, we define
the \emph{proper-walk connection number} $\pW(G)$ as the 
minimum number of colors if one is allowed any
properly colored walk. For symmetry, we will use $\pP(G)$ to
denote the proper connection number.
Trivially, $\pW(G) \le \pP(G)$.

We proceed as follows.
In Section~\ref{s:three} we show that for any connected cyclic graph the proper-walk 
connection number is at most three, and 
in Section~\ref{s:bipartite} we characterize the bipartite graphs that 
have proper-walk and proper connection numbers two. Thereafter,
we show in Section~\ref{s:twoOdd} that the parameter is two
for any graph with two disjoint odd cycles  
and in Section~\ref{s:bridgeless} that the parameter is two for
any bridgeless graph. In Section~\ref{s:forward}
we provide some thoughts on the general case.
We conclude with a comment
about the directed version and some
thoughts for future work.
  
\section{An Upper Bound} \label{s:three}

It is immediate that a properly colored walk cannot 
use the same edge twice in succession. It follows that, in a
tree, every properly colored walk is a path. As  
observed in~\cite{ALLZ16proper}, for the property in trees, one needs the edges 
of the tree to be properly colored, and thus:

\begin{OBS}
If $T$ is  a tree with maximum degree $\Delta$, then
 $\pW(T)=\pP(T)=\Delta$.
 \end{OBS}
 
We present next a general upper bound on the proper-walk connection
number of cyclic graphs. 

\begin{THM}
Let $G$ be a connected graph that is not a tree.
Then $\pW(G) \le 3$.
\end{THM}
\begin{Proof}
We may assume that $G$ is unicyclic (else take suitable spanning subgraph).
Consider the cycle $C$.  Take any proper
coloring of the cycle $C$. For every vertex $v$ of the cycle,
it is incident with two colors in the cycle; so let 
all other edges incident with $v$ have the third color. Color 
the remaining edges so that for every vertex $w$ not on the cycle,
the path $J_w$ from $w$ to the closest vertex of $C$ is properly colored.

There is a properly colored walk between every pair $u$ and $v$ of vertices.
For example, if both $u$ and $v$ are off the cycle, then use $J_u$ to
get to the cycle, go around the cycle to the vertex closest to $v$, and then
use $J_v$ in reverse to get to~$v$.
\end{Proof}

Figure~\ref{f:pw3} gives an example of a graph $G$ where $\pW(G)=3$.
(For a proof of this, see Theorem~\ref{t:cycleFeet}.)

\begin{myFigure}
\includegraphics{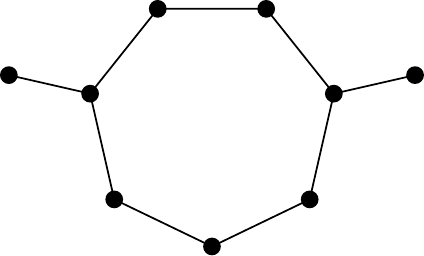}
\\
\caption{A graph $G$ such that $\pW(G)=3$}
\label{f:pw3}
\end{myFigure}

Note that the complete graph has $\pW(G)  = \pP(G)=1$,
while noncomplete graphs have $\pW(G)\ge 2$.
So the big question is: for which graphs is $\pW(G)=2$?

\section{Bipartite Graphs}  \label{s:bipartite}

We next determine which bipartite graphs $G$ have $\pW(G)=2$.

For graph $G$, define $M(G)$ as the spanning subgraph
that results if one removes all the bridges of $G$.
Note that each component of $M(G)$ is either an isolated
vertex or is $2$-edge-connected.

\begin{THM} \label{t:bipartite}
Let $G$ be a connected bipartite graph with order at least $3$. Then $\pW(G)=2$ 
if and only if every component of $M(G)$ is
incident with at most two bridges.
\end{THM}
\begin{Proof}
(1) 
Assume that every component of $M(G)$ is
incident with at most two bridges. We will color the edges of $G$ with two colors.

We first color the edges of $M(G)$.
Let $H$ be a nontrivial component of $M(G)$.
Then $H$ is $2$-edge connected.
By Robbins' Theorem~\cite{Robbins39traffic}, 
such a graph has a strongly connected
orientation, say~$\vec{H}$. (That is, an orientation such that one can get 
from every vertex to every other vertex respecting the orientation.)
Give the vertices of the subgraph $H$ their bipartite coloring; 
then color each arc of~$\vec{H}$ by the color of the head. It follows that 
all directed walks in the orientation~$\vec{H}$ alternate colors.
And within the undirected $H$, each pair of vertices is joined by a properly colored
walk that starts with any designated color or that ends with any
designated color, by either following the arcs or going against the arcs.
Do this process for all nontrivial components $H$ in $M(G)$.

We next color the bridges of $G$. We will color them such that 
two bridges incident with the same component $H$ of $M(G)$ 
have the same color if and only if their ends in $H$ are in different 
partite sets.  This can be achieved by 
considering the graph $F$ obtained from $G$ by contracting each 
(nontrivial) component $H$ of~$M(G)$ to a single vertex $c_H$. 
Note that 
$F$ is acyclic, has maximum degree at most $2$, and is connected, so that 
$F$ is a path. Each edge in $F$ corresponds to a bridge in $G$;
for each bridge $e$ of $G$, let $e'$ be the corresponding
edge in~$F$.

We color $F$ as follows. Start at a leaf-edge
and color it arbitrarily. For subsequent edges, suppose that edge $e'$ is
colored and we need to color adjacent edge $f'$. Say 
edges $e'$ and $f'$ have common end $c_H$ in $F$. Then let $v_e$ be the end of $e$ in $H$ 
and similarly with $v_f$. If $v_e$ and $v_f$ are in the same partite
set of $G$, then give edges $e$ and $f$ different colors; and
if $v_e$ and $v_f$ are in different partite
sets of $G$, then give edges $e'$ and $f'$ the same color.
Finally, transfer the coloring of $F$ to $G$; that is, give each bridge $e$ of $G$ the color of its
corresponding edge $e'$ in $F$.

We claim the above coloring has the desired property; that is, there is a properly colored
walk between every two vertices $u$ and $v$ of $G$.

Every path from $u$ to $v$ uses the same set of bridges in the
same order. So consider two consecutive of these bridges, say 
$b_1$ and $b_2$. Then there is a component of $M(G)$, say $H$,
to which they are both incident.
By the way we colored the bridges, if $b_1$ and $b_2$ have the 
same color, say red, then their ends in~$H$ are in different partite
sets and so every path between them finishes with the same
color it starts with. By above there
is a path between those two ends starting and finishing with a 
blue edge. On the other hand, if $b_1$ and $b_2$ have different
colors, say the former is red and the latter is blue, then these
ends are in the same partite set, and so every path between them finishes
with color different to its start.
By above there is a path between these ends starting with
a blue edge (and necessarily ending with a red edge).
Thus we can piece together the bridges with suitable paths in $M(G)$
to obtain the alternating $u$--$v$-walk

(2) Assume that $G$ has a suitable $2$-coloring.
Since $G$ is bipartite, every closed walk has
the same parity. So assume a properly colored walk enters a subgraph~$H$ of 
$M(G)$ along bridge $b_1$ to vertex~$v_1$ and exits
$H$ along bridge $b_2$ from vertex~$v_2$ (with $v_1=v_2$ allowed). Then
$b_1$ and~$b_2$ must have color determined by
the parity of the distance between $v_1$ and~$v_2$. That is, bridges $b_1$ and $b_2$ have the same color
if and only if $v_1$ and $v_2$ are in different partite sets in $G$.

So suppose there are three bridges $b_1,b_2,b_3$ incident with 
(not necessarily distinct) vertices $v_1,v_2,v_3$ of
$H$. Without loss of generality, $v_1$ and $v_2$ are in the same
partite set~$X$. Thus $b_1$ and $b_2$ need different colors.
But then if $v_3$ is in~$X$, the bridge $b_3$ 
needs a color different from both $b_1$ and $b_2$; and if 
$v_3$ is in the other partite set, then $b_3$ needs to be
the same as both $b_1$ and $b_2$; in each case an impossibility.
\end{Proof}

It turns out that the above characterization also holds
for the proper connection number. 
For,
 in a bipartite graph, all closed walks have even length. Thus, if the edges
are $2$-colored, then there is a properly colored walk between two vertices
if and only if there is a properly colored path between them.  That is:

\begin{THM}
Let $G$ be a connected bipartite graph of order at least $3$. Then $\pP(G)=2$ 
if and only if every component in $M(G)$ is
incident with at most two bridges.
\end{THM}
\noindent
It was known that $\pP(G)=2$ for bridgeless bipartite graphs~\cite{BFGMMMT12proper}.

\section{Disjoint Odd Cycles}  \label{s:twoOdd}

We now consider the general problem of which graphs $G$ have $\pW(G)=2$.

\begin{THM} \label{t:twoOddCycles}
If a connected noncomplete graph $G$ has two edge-disjoint odd cycles, then $\pW(G)=2$.
\end{THM}
\begin{Proof}
Since the graph is noncomplete, we need at least two colors.

Let $C_1$ and $C_2$ be edge-disjoint odd cycles. 
 If they are also
vertex-disjoint, let $P$ be a shortest path joining them; say $P$ starts with vertex $u_1$
in $C_1$ and ends at $u_2$ in $C_2$. If the cycles have a vertex in common, then let
$u_1=u_2$ be such a vertex. Let $H$ be the subgraph consisting of $C_1$, $C_2$, and
$P$ if needed.

Now, color the two edges of $C_1$ incident with~$u_1$ red; then color 
the remaining edges of $C_1$ alternating red and blue so that $u_1$ is
the only vertex not incident with an edge of each color. 
Further, if $P$ exists, color the edges of $P$ alternating colors so that the edge incident with $u_1$ is blue. 
Now, if $P$ has even length or the cycles had a vertex in common, 
color the two edges of $C_2$ incident with~$u_2$ blue; then color 
the remaining edges of $C_2$ alternating red and blue so that $u_2$ is
the only vertex not incident with an edge of each color. 
On the other hand,
if $P$ has odd length, then proceed similarly, except that the two edges of $C_2$
incident with $u_2$ are colored red. 

We claim that this coloring has the property that between every pair $u$ and~$u'$ of (not
necessarily distinct) vertices in $H$, there is a properly colored walk that starts and finishes 
with any prescribed colors. To see this, first note that every vertex of $H$ is incident 
with at least one edge of each color. Thus one can start walking from $u$ with any prescribed color.
Then one can extend this alternating walk indefinitely such that eventually one traverses $P$ in both
directions (if it exists), and goes around both $C_1$ and $C_2$ in both directions. Using this, one can arrive at vertex $u'$ 
having just traversed any designated incident edge.

Now consider the vertices not in $H$.
By choosing a spanning subgraph if needed,
one may assume that for each vertex $v$ not in $H$ there is a unique
path~$J_v$ from $v$ to $H$.
Color the remaining edges such that each $J_v$
is properly colored. See Figure~\ref{f:twoOddCycles}.

\begin{myFigure}
\includegraphics{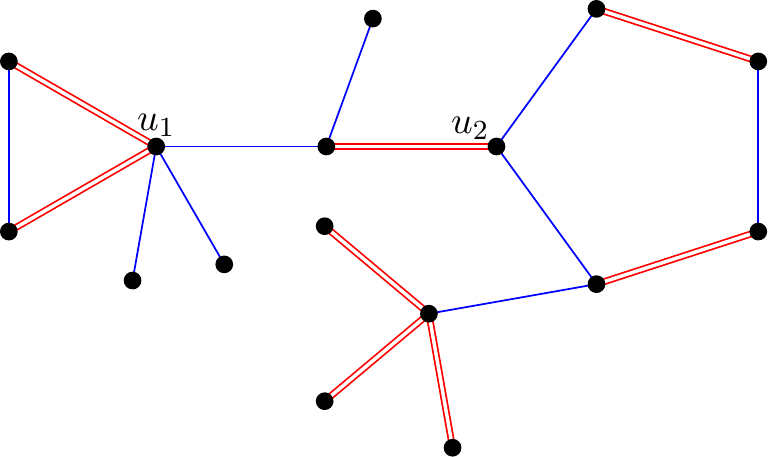}
\\
\caption{A coloring of a graph with two disjoint odd cycles}
\label{f:twoOddCycles}
\end{myFigure}

We claim that the coloring has the desired property. 
To get between any two vertices $v$ and $w$ in $G$, use the alternating path $J_v$
to get to $H$ if necessary, go around $H$ in the appropriate direction, and 
then use the alternating path $J_{w}$ in reverse if needed.
\end{Proof}

Our focus is on simple graphs, but we consider in passing 
what happens if the graph has loops. It is immediate
from the above that if the graph has two loops  
then $\pW(G)=2$, as one can treat the loops as odd cycles.
But actually, $\pW(G)=2$ for any graph with a loop. For, 
one can color the loop blue say, the edges incident with the
loop red, and then alternate colors away from the loop.
There is a properly colored walk between every pair of vertices
by going via the loop.

\section{Bridgeless Graphs}  \label{s:bridgeless}

In this section, we show that $\pW(G) \le 2$ for all connected
graphs $G$ without 
bridges.

\subsection{Preliminaries}

We will need the following simple observation.

\newpage

\begin{OBS} \label{ob:inverterOdd}
Let $P$ be an induced path. (That is, the subgraph induced by the
vertices of $P$ is a path.)
If there is an odd cycle that shares at least one edge with the path $P$,
then there exists a nontrivial path $S$ that is internally disjoint from~$P$
and creates an odd cycle with~$P$.
\end{OBS}
\begin{Proof}
Let $C$ be any odd cycle that shares an edge with $P$.
Consider the vertices of $C\cap P$. Since $P$ is induced,
there must be at least one vertex in $C$ not on $P$. Since
$C$ and $P$ share an edge, there are at least two
vertices in $C\cap P$.
Now, partition the edges of $C$ not in $P$ into segments, where the ends of a segment
are in $P$ and internal vertices of each segment are not in $P$.
If every segment creates an even cycle with $P$, then the result is
bipartite, a contradiction. So some segment creates an odd
cycle with $P$, as required.
\end{Proof}

We will also need the following result.

\begin{LEM} \label{l:orientable}
Let $P$ be a path in a graph $G$ from vertex 
$u$ to vertex $v$ such that for every vertex $w$ not on $P$ there are 
two internally disjoint paths from $w$ to $P$ ending at different vertices.
Then one can orient $G$ such that: \\
(a) $P$ is oriented from $u$ to $v$; \\
(b) one can get from $u$ to every other vertex $w$ by a directed walk; and \\
(c) for all vertices $w_1$ and $w_2$ there is a directed walk between
them in at least one direction.
\end{LEM}
\begin{Proof}
We will create a spanning oriented subgraph $\vec{H}$ such that 
for each vertex $w$ not on $P$:
there exist distinct vertices $q_w$ and $r_w$ on $P$, 
with $q_w$ nearer to $v$, such that there is a directed walk from $q_w$ to $w$
and a directed walk from~$w$ to $r_w$.

Start with $\vec{H}$ as the path $P$ oriented from $u$ to $v$. 
We will grow $\vec{H}$ to contain all the vertices.
Let $w$ be any vertex not on $P$.
Since $P$ contains all the cut-vertices of $G$, there are
two internally disjoint paths from $w$ to $P$.
Say these paths end at vertices $h_1$ and~$h_2$, where $h_1$ is nearer to $v$. 
Add all the vertices of both these paths to $H$, and 
orient the path between $w$ and $h_1$ towards $w$ while orienting the path between $w$ and $h_2$ 
away from $w$. For all newly added vertices, $h_1$ is the $q$-vertex and $h_2$ the $r$-vertex.

If there is still a vertex not in $H$, let $w'$ be such a vertex.
Take the two internally
disjoint paths from~$w'$ to $P$ and cut each when it reaches a vertex that is 
already in $H$. Say we have internally disjoint paths $L_1$ and $L_2$ from~$w'$ to  
vertices $k_1$ and $k_2$. For convenience, if vertex $k$ is on $P$ then we define
$r_k = q_k = k$. 
By reordering $k_1$ and $k_2$ if necessary, it follows that we may assume
$q_{k_1}$ is strictly nearer to~$v$ than~$r_{k_2}$. 

Add all the vertices of both paths $L_1$ and $L_2$ to $H$.
Then orient $L_1$ towards $w'$ and orient $L_2$ away from $w'$. 
For all newly added vertices, $q_{k_1}$ is the $q$-vertex and $r_{k_2}$ the $r$-vertex.
Repeat this procedure until $H$ contains all the vertices.

We claim this orientation $\vec{H}$ has the desired three properties. The first property was explicitly
satisfied. For a directed walk from $u$ to $w$, go along $P$ to~$q_w$ and then along the 
walk to $w$.  
Further, without loss of generality, we may assume that 
$r_{w_2}$ is not farther from $v$ than $r_{w_1}$;
this means that $q_{w_2}$ is nearer to~$v$ than~$r_{w_1}$. So one can 
get from $w_1$ to $w_2$ by going to $r_{w_1}$, going along $P$ to $q_{w_2}$,
and thence to $w_2$.
\end{Proof}

\subsection{Main Result}

We define a \dd{theta-graph} as a graph that is formed by taking
a cycle $C$ of even length (called the outer cycle) and a path 
$P$ (called the inverter) and identifying the ends of the path $P$
with two vertices $u$ and $v$ of the cycle $C$ such that the
result is nonbipartite. See Figure~\ref{f:theta} for an example.

\begin{myFigure}
\includegraphics{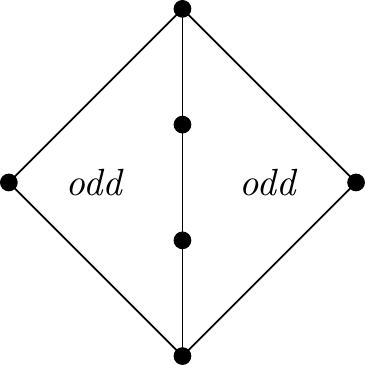} 
\\
\caption{A theta-graph}
\label{f:theta}
\end{myFigure}

\newpage

\begin{THM}
There does not exist a $2$-connected graph $G$ such that $\pW(G)=3$.
\end{THM}
\begin{Proof}
Suppose block $G$ has $\pW(G)=3$.
We saw above (Theorem~\ref{t:bipartite}) that $G$ cannot be bipartite.
Also, we saw (Theorem~\ref{t:twoOddCycles}) that $G$ does not contain two
edge-disjoint odd cycles.

Consider some odd cycle of the graph $G$. If it is a hamilton cycle, 
then it is easily seen that $\pW(G)=2$. So assume there is a vertex not on this cycle.
By $2$-connectedness, we can find two disjoint paths from this vertex
to the cycle, ending at vertices $u$ and $v$ say. That is, we have three internally
disjoint $u$--$v$ paths such that the result is not bipartite. Two of 
these paths have the same parity; choose them to be the outer cycle,
and the other path to be the inverter. That is, the result is a theta-graph. 

Out of all theta-subgraphs, 
\begin{quote} \slshape
choose the theta-subgraph where the inverter $P$ is as short as 
possible.
\end{quote}
Let $C$ be the outer cycle of the chosen theta-graph.

\begin{claim} \label{c:isBipartite}
The graph $G-C$ is bipartite. 
\end{claim}
\begin{Proof}
Suppose there is an odd cycle in $G-C$. Since there are not
two edge-disjoint odd cycles in $G$, that odd cycle must share an edge with (the interior of) 
$P$. Then by Observation~\ref{ob:inverterOdd}, there is segment $S$ in $G-(C\cup P)$ 
that joins two vertices of the interior of $P$ but is otherwise disjoint from $P$ and
creates an odd cycle with $P$. This segment $S$ combined
with $P$ and either half of $C$ provides a theta-graph with a shorter
inverter, which contradicts our choice of theta-subgraph.
\end{Proof}

Let $P'$
be the path $P$ minus $u$ and $v$.
Partition the vertices not in the theta-graph into two sets: let $A$
be those vertices that can 
reach the outer cycle $C$ without going through $P'$, 
and let~$B$ be those that cannot.

Now, color the graph $G$ as follows. Color the theta-graph such that 
the outer cycle $C$ is properly colored, as is the inverter $P$. 
Without loss of generality, assume that $C$ is drawn so that 
every properly colored walk leaving the inverter proceeds clockwise on the outer
cycle.

For each vertex $w$ of $A$, retain one path $J_w$ to $C$ that does not intersect $P'$. 
Color the edge of $J_w$ incident with $C$ 
such that one can go across that edge and proceed counter-clockwise
around the outer cycle.
Color the remaining edges of the path $J_w$ so that it is properly colored. 

Finally, consider the set $B$. 
By Claim~\ref{c:isBipartite}, the graph  $H$ induced by $P' \cup B$ is bipartite.
Since the graph $G$ is $2$-connected, there are two internally disjoint
paths from every vertex $w\in B$ to the theta-graph. By the definition of $B$,
these paths must meet the theta-graph on $P'$. Thus we can apply Lemma~\ref{l:orientable}
to $H$ and~$P'$ to obtain an orientation $\vec{H}$ with the properties listed in that lemma. 
Give each vertex of $H$ its bipartite coloring; then color each arc of 
$\vec{H}$ by the color of its head.
As we used in the proof of Theorem~\ref{t:bipartite}, 
in such a coloring every walk that respects the orientation automatically alternates
colors. 

See Figure~\ref{f:theta-graph} for an example, where
the vertices of $A$ are drawn outside the outer cycle and the vertices of~$B$ are
drawn inside the outer cycle.

\begin{myFigure}
\includegraphics{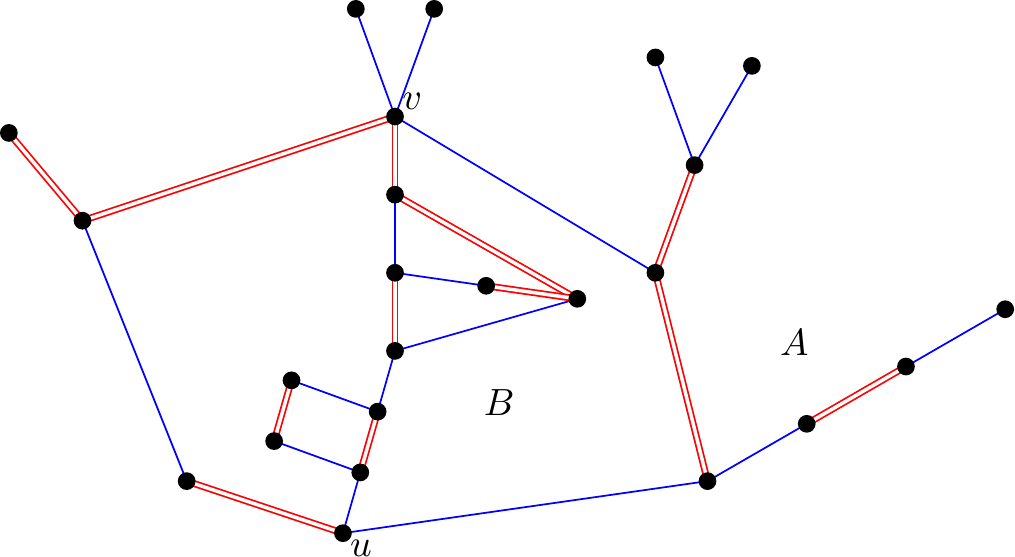}
\\
\caption{Coloring of theta-graph and a spanning subgraph} 
\label{f:theta-graph}
\end{myFigure}

We claim that the resultant coloring of $G$ has a properly colored walk between any pair of 
vertices. For example, to get from a vertex $w_1$ of $A$ to another vertex $w_2$ of $A$,
follow $J_{w_1}$, go counter-clockwise around the outer cycle, over the inverter, clockwise
around the outer cycle, and then use $J_{w_2}$ in reverse. 
To get from a vertex $w_1$ of $A$ to a vertex $w_3$ of $B$, 
follow $J_{w_1}$, go counter-clockwise around the outer cycle to $u$, and then 
over the inverter to the first vertex of $P'$, and then use the directed walk in $\vec{H}$
to $w_3$. And, to get between vertices $w_3$ and $w_4$ of $B$, use the 
directed walk in $\vec{H}$.
\end{Proof}

From the above result, the question of bridgeless graphs is easily resolved:

\begin{THM} \label{t:bridgeless}
If $G$ is a connected bridgeless graph, then $\pW(G)\le 2$.
\end{THM}
\begin{Proof}
Assume $G$ is bridgeless but not $2$-connected. Consider the blocks
of~$G$. If any two of these are nonbipartite, then there
are two edge-disjoint odd cycles, and the result follows 
from Theorem~\ref{t:twoOddCycles}. 
If all the blocks are bipartite, then
the result follows from Theorem~\ref{t:bipartite}. 
So assume that exactly one block, say~$H$, is not bipartite.

By the above theorem, that block $H$ can be colored with two colors to have 
a properly colored walk between every pair of vertices in $H$. 
Color all remaining blocks properly, as in Theorem~\ref{t:bipartite}.
We claim the resultant coloring has the desired property. To find a properly
colored walk between vertices $u$ and $v$, let $u'$ be the vertex of $H$
nearest to $u$ and $v'$ the vertex of $H$ nearest to $v$. Then find
the properly colored walk between $u'$ and $v'$. This can be extended
to a properly colored walk between $u$ and $v$, since there is a walk from $u$ to $u'$
ending with any prescribed color, and a walk from $v'$ to $v$ starting with
any desired color.
\end{Proof}

\section{Unicyclic Graphs} \label{s:forward}

It is unclear what happens in general in graphs with bridges. The precise placement
of bridges seems to matter. For example, consider
the collection $\calG$ of graphs formed by taking an odd cycle and adding
feet to some of the vertices of the cycle. (By adding a foot we mean
adding a new vertex and joining it to exactly one vertex of the cycle.)

\begin{THM} \label{t:cycleFeet}
Let $G$ be a graph of $\calG$. Then $\pW(G)=2$ if and only 
if there are three consecutive vertices $u,v,w$ on the cycle
such that $u$ is adjacent to at most one foot, $w$ is adjacent to at most
one foot, and all vertices other than $u,v,w$ are  incident with no feet.
\end{THM}

\begin{myFigure}
\begin{center}
\includegraphics{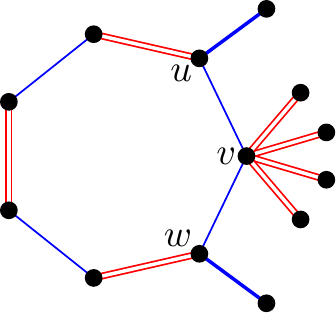} 
\\
\caption{A graph in $\calG$ with $\pW=2$}
\label{f:calGPW2}
\end{center}
\end{myFigure}

\begin{Proof}
(1) We first prove that the conditions are necessary for the graph to have $\pW=2$.
That is, assume the graph has a $2$-coloring
such that every vertex can reach every other vertex by a properly colored
walk.

\textit{Case~A: Assume the odd cycle has length at least $5$.}\\
Call a vertex of the cycle a break if the two cycle edges incident with it have the same
color. The number of breaks has the same parity as the number of vertices;
that is, there is an odd number of breaks. Suppose there are
at least three breaks, say vertices $v_1$, $v_2$, and $v_3$. These divide
up the cycle into three paths, at least one of which must have more than one edge,
say the $v_1$--$v_2$ path. Then there is no alternating walk between $v_3$ and an interior
vertex of that path, a contradiction. That is, there is exactly one break.

So let $v$ be the unique vertex on the cycle incident
with two edges of the same color, with neighbors $u$ and $w$ on the cycle. 
Suppose there is a foot attached to a vertex $x$
that is neither $u$, $v$, nor $w$. Then the two edges of
the cycle incident with $x$ have different colors, and so any 
walk from the foot can proceed in only one direction around the 
cycle, and gets stopped at $v$ without reaching all the vertices.
Thus, all feet must be attached to one of $u$, $v$, or $w$. 

Consider a foot incident with $u$. In order for it to reach all
vertices, the edge incident with it must have the same color as the $uv$ edge.
It follows that the foot is unique, since otherwise the two feet would not be able
to reach each other.

\textit{Case~B: Assume the odd cycle is a triangle.}\\
If the triangle has
exactly one break, then 
by the same argument as Case~A, the other two vertices of the triangle can 
be incident with at most one foot each. Further, if the triangle is monochromatic,
then it is easy to see that each vertex of the cycle is incident with at most one foot.
\smallskip

(2) We second prove that the conditions are sufficient.
Color the cycle such that $v$ is incident with two edges of the same color
and every other vertex sees both colors. Color the leaf incident
with $u$ and/or $w$ with the same color as the $uv$ edge;
color all leaves incident with $v$ with the other color.
It is easily checked that this coloring has the desired property.
\end{Proof}

\section{Directed Graphs} \label{s:digraph}

For a strongly connected digraph, one can define the proper-walk
connection number as in the undirected case. 
 This idea was recently
introduced for paths by Magnant et al.~\cite{MMPSW16directed}.
They showed that:

\begin{THM} \cite{MMPSW16directed}
Let $D$ be a strongly-connected digraph.
Then $\pP(D) \le 3$.
\end{THM}

This is sharp, even for the proper-walk case, since 
an odd cycle needs three colors; that is, $\pW(D)= 
\pP(D) = 3$ if $D$~is an odd cycle.

We note that the two parameters can be different. That is, 
there are digraphs 
with $\pW(D)=2$ and $\pP(D)=3$.
For example, take two disjoint directed triangles and identify one vertex
of each. See Figure~\ref{f:bowTie}.

\begin{myFigure}
\includegraphics{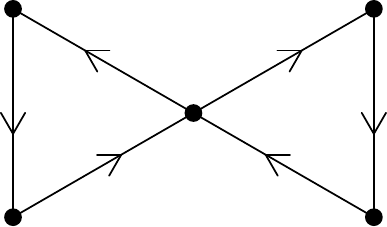}
\\
\caption{A graph where $\pW( D) = 2$ and $\pP(D) = 3 $}
\label{f:bowTie}
\end{myFigure}

Another direction is to add loops. If one adds loops
at all vertices, then one needs only two colors  (color all original
arcs one color and all loops a second color).

\section{Conclusion} \label{s:conclusion}

We proved that every connected graph has proper-walk connection number
at most three, and showed that it is two for some families.
One natural 
open problem is the complexity of recognizing which
graphs have the parameter $2$. Is there a polynomial-time algorithm,
or is it NP-hard? Note that it is easy to check using a breadth-first-search 
whether a given coloring
has a properly colored walk  between two vertices.

Other directions of interest include the question where some of the
edges of the graph are already colored. For example, K\'{e}zdy and Wang 
\cite{KW99alternating} asked when one could complete a $2$-coloring
such that there is an alternating path between two specified vertices.
One could also insist on stronger properties; for example, that every pair 
of vertices is in a properly colored cycle, 
or closed walk.

\section*{Acknowledgements}

We thanks the referees for their thoughtful comments that improved the paper.


\begin{thebibliography}{9}
\bibitem{ALLZ16proper}
E.~Andrews, C.~Lumduanhom, E.~Laforge, P.~Zhang, 
On proper-path colorings in  graphs, 
\textit{J. Combin. Math. Combin. Comput.} 97 (2016), 189--207.

\bibitem{BFGMMMT12proper}
V.~Borozan, S.~Fujita, A.~Gerek, C.~Magnant, Y.~Manoussakis, L.~Montero, Z.~Tuza, 
Proper connection of graphs, 
\textit{Discrete Math.} 312 (2012), 2550--2560.

\bibitem{GLQ16random}
R.~Gu, X.~Li, Z.~Qin, 
Proper connection number of random graphs, 
\textit{Theoret. Comput. Sci.} 609 (2016), 336--343.

\bibitem{KW99alternating}
A.~K{\'e}zdy, C.~Wang, 
Alternating walks in partially {$2$}-edge-colored graphs and optimal strength of graph labeling, 
\textit{Discrete Math.} 194 (1999), 261--265.

\bibitem{LM15survey}
X.~Li, C.~Magnant, Properly colored notions of connectivity---a dynamic survey,
\textit{Theory and Applications of Graphs} 0 (2015), Article 2.

\bibitem{MMPSW16directed}
C.~Magnant, P.~Morley, S.~Porter, P.~Salehi~Nowbandegani, H.~Wang, 
Directed proper connection of graphs, 
\textit{Mat. Vesnik} 68 (2016), 58--65.

\bibitem{Robbins39traffic}
H.~E. Robbins, A theorem on graphs, with an application to a problem of traffic control, 
\textit{Amer. Math. Monthly} 46 (1939), 281--283.

\end{thebibliography}
\end{document}